%% file: MT-HVDC_AC_Control_secondary_converter_conference_IFAC_draft4.tex
\documentclass{ifacconf}
\usepackage{IEEEtrantools}
\usepackage{natbib} 
\usepackage[english]{babel}

\usepackage{mathrsfs}
\usepackage{amsfonts}
\usepackage{amssymb}
\usepackage{amsmath}
\usepackage{graphicx, caption}
\usepackage{commath}
\usepackage{tikz}
\usepackage{tkz-graph}
\usepackage[europeanresistors]{circuitikz}
\usetikzlibrary{arrows}
\usepackage{flushend}
\usepackage{psfrag}
\usepackage{pgfplots}
\usepackage{todonotes}
\usetikzlibrary{external}
\usepackage{blockdiagram}
\usepackage{caption}
\usepackage{subcaption}
\usepackage{color}
\usepackage{booktabs}

\usetikzlibrary{external}
\pgfplotsset{every x tick label/.append style={font=\small, yshift=0.2ex}} 
\pgfplotsset{every y tick label/.append style={font=\small, xshift=0.3ex}} 
\pgfplotscreateplotcyclelist{Necsys}{%
	{blue,line width=0.5pt},
	{red,line width=0.5pt},
	{black!30!green,line width=0.5pt},
	{cyan,line width=0.5pt},
	{yellow!90!black,line width=0.5pt},		
	{magenta,line width=0.5pt},
}

\newtheorem{theorem}{Theorem}
\newtheorem{corollary}[theorem]{Corollary}

\newtheorem{remark}{Remark}

\newtheorem{assumption}{Assumption}
\newtheorem{objective}{Objective}

\DeclareMathOperator*{\argmin}{argmin}

\DeclareMathOperator*{\diag}{diag}

\newcommand{\beq}{\begin{equation}}
\newcommand{\eeq}{\end{equation}}
\newcommand{\bq}{\begin{eqnarray}}
\newcommand{\eq}{\end{eqnarray}}
\newcommand{\bqn}{\begin{eqnarray*}}
\newcommand{\eqn}{\end{eqnarray*}}
\newcommand{\bee}{\begin{enumerate}}
\newcommand{\eee}{\end{enumerate}}

\begingroup
\makeatletter
\renewcommand{\p@subfigure}{}
\makeatother

\newlength\fheight
\newlength\fwidth

\begin{document}
\begin{frontmatter}

\title{Coordinated Frequency Control through \\ MTDC Transmission Systems \thanksref{footnoteinfo}}

\author[First]{Martin Andreasson} 
\author[Second]{Roger Wiget} 
\author[First]{Dimos V. Dimarogonas}
\author[First]{Karl H. Johansson}
\author[Second]{G\"oran Andersson} 

\address[First]{ACCESS Linnaeus Centre,
KTH Royal Inst. of Technology, Sweden.\\
(e-mail: {\tt \small \{mandreas, dimos, kallej\}@kth.se}.) }

\address[Second]{Power Systems Laboratory, ETH Zurich, Switzerland.\\
(e-mail: {\tt \small \{wiget, andersson\}@eeh.ee.ethz.ch }.) }

\thanks[footnoteinfo]{This work was supported in part by the European Commission, the Swedish Research Council (VR) and the Knut and Alice Wallenberg Foundation.
We would like to thank the anonymous reviewers for their valuable comments and suggestions. 
 Corresponding author: Martin Andreasson, e-mail: mandreas@kth.se.}
 

\begin{abstract}: 
In this paper we propose a distributed dynamic controller for sharing frequency control reserves of asynchronous AC systems connected through a multi-terminal HVDC (MTDC) grid. We derive sufficient stability conditions, which guarantee that the frequencies of the AC systems converge to the nominal frequency. Simultaneously, the global quadratic cost of power generation is minimized, resulting in an optimal distribution of generation control reserves. The proposed controller also regulates the voltages of the MTDC grid, asymptotically minimizing a quadratic cost function of the deviations from the nominal voltages. 
The proposed controller is tested on a high-order dynamic model of a power system consisting of asynchronous AC grids, modelled as IEEE 14 bus networks, connected through a six-terminal HVDC grid. The performance of the controller is successfully evaluated through simulation. 
\end{abstract}

\end{frontmatter}

\section{Introduction}
Power transmission over long distances with low losses is one of the main challenges in today's power transmission systems. As the share of renewables rises, so does the need to balance generation and consumption mismatches, often over large geographical areas.  Due to the high resistive losses in AC cables, high-voltage direct current (HVDC) power transmission is a commonly used technology for long-distance power transmission. 
The higher investment cost of an HVDC transmission system compared to an AC transmission system is compensated by the lower resistive losses for sufficiently long distances \citep{melhem2013electricity}. The break-even point, i.e., when the total construction and operation costs of overhead HVDC and AC lines equal, is typically 500--800 km \citep{padiyar1990hvdc}. However, for cables, the break-even point is typically less than 50 km \citep{Hertem2010technical}. Increased use of HVDC technologies for electrical power transmission suggests that future HVDC transmission systems are likely to consist of multiple terminals connected by several HVDC transmission lines \citep{Haileselassie2013Power}. Such systems are referred to as multi-terminal HVDC (MTDC) systems in the literature.
Many existing AC grids are connected through legacy HVDC links, which are typically used for bulk power transfer between AC areas, rather than balancing generation mismatches. The fast operation of the DC converters however also enables frequency regulation of one of the connected AC grids through the HVDC link. One example is the frequency regulation of the island of Gotland in Sweden, which is connected to the Nordic grid through an HVDC cable \citep{axelsson2001gotland}. However, since the Nordic grid has orders of magnitudes higher inertia than the AC grid of Gotland, the influence of the frequency regulation on the Nordic grid is negligible.

By connecting multiple AC grids by an MTDC system, the frequency regulation reserves in each AC grid can be shared with connected AC grids, which reduces the need for frequency regulation reserves in the individual AC systems \citep{li2008frequency}.  
In recent years, this idea has gained increasing interest in the literature. 
 \cite{dai2011voltage} and \cite{silva2012provision} employ several controllers with decentralized structure to share frequency control reserves. In \citep{silva2012provision} no stability analysis of the closed-loop system is performed, whereas \citep{dai2011voltage} guarantees stability provided that the connected AC areas have identical parameters and the voltage dynamics of the HVDC system are neglected. In \citep{taylordecentralized}, an optimal decentralized controller for AC systems connected by HVDC systems is derived. In contrast to the aforementioned references, \cite{andreasson2014distributed} consider the dynamics of connected AC systems as well as the dynamics of the MTDC system.
By connecting the AC areas with a communication network, the performance of the controller can be further improved. 
\cite{dai2010impact} consider a distributed controller relying on a communication network to share frequency control reserves of asynchronous AC transmission systems connected through an MTDC system. However, the controller requires a slack bus to control the DC voltage, and is thus only able to share the generation reserves of the non-slack AC areas. 
Another distributed controller is proposed by \cite{dai2013voltage}. Stability is guaranteed, and the need for a slack bus is eliminated. The voltage dynamics of the MTDC system are however neglected. Moreover the implementation of the controller requires every controller to access measurements of the DC voltages of all other MTDC terminals. 
 \citep{andreasson2014distributedSecondary} propose a distributed secondary generation controller. In contrast to the aforementioned references, the MTDC dynamics are explicitly modeled, and the DC voltages are controlled in addition to the frequencies. 

Despite the coordination through a communication network, the distributed controller by \cite{andreasson2014distributedSecondary} fails in eliminating the static errors of the frequencies. In this paper we address this issue by introducing a distributed secondary controller also for the power fed into the MTDC grid from the AC grids. This allows us to eliminate any static errors in the frequency deviations. Furthermore, quadratic cost functions of the power generation and the voltage deviations are minimized. 

 The remainder of this paper is organized as follows. In Section \ref{sec:model}, the system model and control objectives are defined. In Section \ref{sec:secondary_frequency_control}, a distributed secondary frequency controller for sharing frequency control and restoration reserves is presented, and is shown to satisfy the control objectives. 
 In Section \ref{sec:simulations}, simulations of the controller on a six-terminal MTDC test system are provided, showing the effectiveness of the proposed controller. The paper ends with concluding remarks in Section~\ref{sec:discussion}.


\section{Model and problem setup}
\label{sec:model}
\subsection{Notation}
\label{subsec:prel}
Let $\mathcal{G}$ be a static, undirected graph. Denote by $\mathcal{V}=\{ 1,\hdots, n \}$ the vertex set of $\mathcal{G}$, and by $\mathcal{E}=\{ 1,\hdots, m \}$ the edge set of $\mathcal{G}$. Let $\mathcal{N}_i$ be the set of neighboring vertices to $i \in \mathcal{V}$.
Denote by $\mathcal{B}$ the vertex-edge incidence matrix of $\mathcal{G}$, and let $\mathcal{\mathcal{L}_W}=\mathcal{B}W\mathcal{B}^T$ be the weighted Laplacian matrix of $\mathcal{G}$, with edge-weights given by the  elements of the diagonal matrix $W$. We denote the space of real-valued $n\times m$-valued matrices by $\mathbb{R}^{n\times m}$.
Let $\mathbb{C}^-$ denote the open left half complex plane, and $\bar{\mathbb{C}}^-$ its closure. We denote by $c_{n\times m}$ a vector or matrix of dimension $n\times m$, whose elements are all equal to $c$. For simplifying notation, we write $c_n$ for $c_{n\times 1}$. For a symmetric matrix $A$, $A>0 \;(A\ge 0)$ is used to denote that $A$ is positive (semi) definite. $I_{n}$ denotes the identity matrix of dimension $n$. For simplicity, we will often drop the notion of time dependence of variables, i.e., $x(t)$ will be denoted $x$. Let $\norm{\cdot}_\infty$ denote the maximal absolute value of the elements of a vector.
\subsection{Model}
\label{subsec:model}
We will give here a unified model for an MTDC system interconnected with several mutually asynchronous AC systems.
We consider an MTDC transmission system consisting of $n$ converters, denoted $i=1, \dots, n$, each connected to an AC system, i.e., there are no pure DC nodes of the MTDC grid. The converters are assumed to be connected by an MTDC transmission grid. The dynamics of converter $i$ is assumed to be given by
\begin{align}
\begin{aligned}
C_i \dot{V}_i &= -\sum_{j\in \mathcal{N}_i} \frac{1}{R_{ij}}(V_i -V_j) + I_i^{\text{inj}} ,
\end{aligned}
\label{eq:hvdc_coordinated_voltage}
\end{align}
where $V_i$ is the voltage of converter $i$, $C_i>0$ is its capacitance, and $I_i^{\text{inj}}$ is the injected current from an AC grid connected to the DC converter.  The constant $R_{ij}$ denotes the resistance of the HVDC transmission line connecting the converters $i$ and $j$. 
The graph corresponding to the HVDC line connections is assumed to be connected.
 The AC system is assumed to consist of a single generator which is connected to the corresponding DC converter, representing an aggregate model of the AC grid. The dynamics of the AC system are given by the swing equation \citep{machowski2008power}:
\begin{align}
m_i \dot{\omega}_i &=  P^\text{gen}_i  + P_i^{{m}} - P_i^{\text{inj}}, \label{eq:frequency}
\end{align}
where $m_i>0$ is its moment of inertia. The constant $P^\text{gen}_i$ is the generated power, $P^m_i$ is the power load  and $P_i^{\text{inj}}$ is the power injected to the DC system through converter $i$, respectively. The powers are all assumed to be deviations from a nominal operation point. 
The control objective can now be stated as follows. 
\begin{objective}
\label{obj:1_hvdc_coordinated}
The frequency deviations are asymptotically equal to zero, i.e.,
\begin{align}
\lim_{t\rightarrow \infty} \omega_i(t)-\omega^{\text{ref}} = 0  \quad  i = 1, \dots, n, \label{eq:hvdc_coordinated_frequency_objective}
\end{align}
where $\omega^{\text{ref}}$ is the nominal frequency. 
The total quadratic cost of the power generation is minimized asymptotically, i.e., $\lim_{t\rightarrow \infty} P_i^\text{gen} = P_i^{\text{gen}*},  \forall i=1, \dots, n$, where
\begin{align}
[P_1^{\text{gen}*}, \dots, P_n^{\text{gen}*}]  =  \argmin_{P_1, \dots, P_n} \frac 12 \sum_{i=1}^n f^P_i P_i^2 \label{eq:hvdc_coordinated_generation_objective} 
\end{align}
subject to $P^\text{gen}_i + P_i^{{m}} - P_i^{\text{inj}} = 0, \; \forall i=1, \dots, n$ and $\sum_{i=1}^n P_i^\text{inj} = 0$, i.e., power balance both in the AC grids and in the MTDC grid, in the absence of power losses. The positive constants  $f^P_i$ represent the local cost of generating power. 
Finally, the voltages are such that the a quadratic cost function of the voltage deviations is minimized asymptotically, i.e., $\lim_{t\rightarrow \infty} V_i = V_i^*, \forall i=1, \dots, n$, where
\begin{align}
 [V_1^*, \dots,  V_n^*] = \argmin_{V_1, \dots, V_n} \frac 12 \sum_{i=1}^n f^V_i (V_i - V_i^\text{ref})^2
 \label{eq:hvdc_coordinated_voltage_objective}
\end{align}
subject to \eqref{eq:hvdc_coordinated_frequency_objective}--\eqref{eq:hvdc_coordinated_generation_objective}, and where the $f^V_i$ is a positive constant reflecting the local cost of voltage deviations, and $V_i^\text{ref}$ is the nominal voltage of converter $i$. 
\end{objective}
\begin{remark}
The minimization of \eqref{eq:hvdc_coordinated_generation_objective} is equivalent to power sharing, where the generated power of AC area $i$ is asymptotically inverse proportional to the cost $f_i^P$. The cost $f_i^P$ can be chosen to reflect the available generation capacity of area $i$.
\end{remark}
\begin{remark}
It is in general not possible that $\lim_{t\rightarrow \infty} V_i(t) = V_i^\text{ref} \; \forall i=1, \dots, n$, since this does not allow for the  currents between the HVDC converters to change by \eqref{eq:hvdc_coordinated_voltage}. Note that the optimal solution to \eqref{eq:hvdc_coordinated_generation_objective} fixes the relative voltages, leaving only the ground voltage as a decision variable of \eqref{eq:hvdc_coordinated_voltage_objective}.
Note also that the reference voltages $V_i^\text{ref}, \; i=1, \dots, n$, are generally not uniform, as is the reference frequency $\omega^\text{ref}$. 
 \end{remark}
 \begin{remark}
 Note that Objective \ref{obj:1_hvdc_coordinated} does not include constraints of, e.g., generation and line capacities. This requires that the perturbations from the operating point are sufficiently small, to guarantee that these constraints are not violated. In case of a large disturbance, e.g., a fault, the constraints might be violated. 
 \end{remark}

\section{Coordinated secondary frequency control}
\label{sec:secondary_frequency_control}

\subsection{Controller structure}
\label{subsec:controlelr_structure}
In this section we propose a distributed secondary frequency controller. In addition to the generation controller proposed by \cite{andreasson2014distributedSecondary}, we also propose a secondary controller for the voltage injections into the HVDC grid. 
The distributed generation controller of the AC systems, which was first given by  \cite{andreasson2014distributedSecondary}, is given by
\begin{align}
P^\text{gen}_i &=- K_i^{\text{droop}} (\omega_i-\omega^{\text{ref}}) -  \frac{K^V_i}{K^\omega_i} K^\text{droop, I}_i \eta_i \nonumber \\
\dot{\eta}_i &= K_i^{\text{droop,I}}(\omega_i-\omega^{\text{ref}}) - \sum_{j\in \mathcal{N}_i} c^\eta_{ij} (\eta_i-\eta_j),
\label{eq:hvdc_coordinated_droop_control_secondary_distributed}
\end{align}
where $K_i^\text{droop}$ and $K_i^\text{droop, I}$ are positive controller parameters. Moreover $c^\eta_{ij} = c^\eta_{ji}>0$, i.e., the communication graph is undirected.
The above controller can be interpreted as a distributed PI-controller, with a distributed consensus filter acting on the integral states $\eta_i$. 
The proposed converter controllers governing the incremental power injections from the AC systems into the HVDC grid are given by
\begin{align}
\label{eq:voltage_control_secondary}
P_i^{\text{inj}} &= K_i^{{\omega}} (\omega_i - \omega^{\text{ref}}) + K_i^{{V}}(V_i^{\text{ref}}-V_i) \nonumber \\ 
& \;\;\;\;+ \sum_{j\in \mathcal{N}_i} c^\phi_{ij} (\phi_i - \phi_j) \nonumber \\
\dot{\phi}_i &= \frac{K^\omega_i}{K^V_i} \omega_i - \gamma \phi_i, 
\end{align}
where  $K_i^V$ and $K_i^\omega$ are positive controller parameters, $\gamma>0$ and $c^\phi_{ij} = c^\phi_{ji}>0$. 
The converter controller \eqref{eq:voltage_control_secondary} can be interpreted as an emulation of an AC network between the isolated AC areas. The auxiliary controller variables $\phi_i$ can be thought of as the phase angles of AC area $i$, governing the power transfer between the areas, if these were connected by AC lines rather than an MTDC grid. In contrast to an AC system, however, the power is fed into the HVDC grid instead of being directly transferred to the connected AC areas. 

The HVDC converter is assumed to be perfect and instantaneous, i.e., injected power on the AC side is immediately and losslessly converted to DC power. Furthermore, the dynamics of the converter are ignored, implying that the converter is assumed to track the output of controller \eqref{eq:voltage_control_secondary} instantaneously. This assumption is reasonable due to the dynamics of the converter typically being orders of magnitudes faster than the primary frequency control dynamics of the AC system \citep{kundur1994power}.
The relation between the injected HVDC current and the injected AC power is thus given by
\begin{align}
V_iI_i^{\text{inj}} = P_i^{\text{inj}}. \label{eq:power-current_nonlinear}
\end{align}
By assuming $V_i=V^{\text{nom}}\; \forall i=1, \dots, n$, where $V^{\text{nom}}$ is a global nominal voltage, we obtain
\begin{align}
V^{\text{nom}}I_i^{\text{inj}} = P_i^{\text{inj}}. \label{eq:power-current}
\end{align}

\subsection{Stability analysis}
\label{sec:secondary_frequency_control_stability}
We are now ready to analyze the stability of the closed-loop system. Define the stacked state vectors as $\hat{\omega}=\omega - \omega^\text{ref}1_n$ and $\hat{V}=V - V^\text{ref}$, where $\omega=[\omega_1, \dots, \omega
_n]^T$, $V=[V_1, \dots, V_n]^T$, $V^\text{ref} = [V_1^\text{ref}, \dots, V_n^\text{ref}]^T$, $\eta = [\eta_1, \dots, \eta_n]^T$ and $\phi = [\phi_1, \dots, \phi_n]$. 
Combining the voltage dynamics \eqref{eq:hvdc_coordinated_voltage}, the frequency dynamics \eqref{eq:frequency} and the generation control \eqref{eq:hvdc_coordinated_droop_control_secondary_distributed}, the converter controller \eqref{eq:voltage_control_secondary} and the power-current relationship \eqref{eq:power-current},  we obtain the closed-loop dynamics
\begin{IEEEeqnarray}{lcl}
\dot{\hat{\omega}} &=& M \Big(- (K^{\text{droop}} + K^\omega) \hat{\omega} + K^V \hat{V} \nonumber \\
&& - {K^V}(K^\omega)^{-1} K^\text{droop, I} \eta - \mathcal{L}_\phi \phi  + P^{{m}}  \Big) \nonumber \\
\dot{\hat{V}} &=& 
\frac{1}{V^{\text{nom}}}E {K}^\omega \hat{\omega} -E\left(\mathcal{L}_R + \frac{K^V}{V^{\text{nom}}} \right) \hat{V} + \frac{1}{V^\text{nom}} E \mathcal{L}_\phi  \phi \nonumber \\
\dot{\eta} &=&  K^{\text{droop,I}}\hat{\omega} - \mathcal{L}_\eta \eta \nonumber \\
\dot{\phi} &=&  (K^V)^{-1}{K}^\omega \hat{\omega}  -\gamma \phi, \label{eq:cl_dynamics_vec_delta_int_coordinated}
\end{IEEEeqnarray}
where 
$M=\diag({m_1}^{-1}, \hdots , {m_n}^{-1})$ is a matrix of inverse generator inertia, 
 $E=\diag(C_1^{-1}, \dots, C_n^{-1})$ is a matrix of electrical elastances,  $\mathcal{L}_R$ is the weighted Laplacian matrix  of the MTDC grid with edge-weights $1/R_{ij}$,  $\mathcal{L}_\eta$ and $\mathcal{L}_\phi$ are the weighted Laplacian matrices of the communication graphs with edge-weights $c^\eta_{ij}$ and $c^\phi_{ij}$, respectively, and $P^m=[P^m_1, \dots, P^m_n]^T$. We define the diagonal matrices of the controller gains as $K^\omega = \diag(K^\omega_1, \dots, K^\omega_n)$, etc. 

Let $y=[\hat{\omega}^T, \hat{V}^T]^T$ define the output of \eqref{eq:cl_dynamics_vec_delta_int_coordinated}. Clearly the linear combination $1_n^T\phi$ is unobservable and marginally stable with respect to the dynamics \eqref{eq:cl_dynamics_vec_delta_int_coordinated}, as it lies in the nullspace of $\mathcal{L}_\phi$. In order to facilitate the stability analysis, we will perform a state-transformation to this unobservable mode. Consider the following state-transformation:
\begin{align}
\phi' = \begin{bmatrix}
\frac{1}{\sqrt{n}} 1_n^T \\ S^T
\end{bmatrix} \phi
\qquad
\phi = \begin{bmatrix}
\frac{1}{\sqrt{n}} 1_n & S
\end{bmatrix} \phi'
\label{eq:transformation_phi}
\end{align}
where $S$ is an $n\times(n-1)$ matrix such that $\left[\frac{1}{\sqrt{n}}1_n \; S\right]$ is orthonormal. By applying the state-transformation \eqref{eq:transformation_phi} to \eqref{eq:cl_dynamics_vec_delta_int_coordinated}, we obtain dynamics where it can be shown that 
the state $\phi_1'$ is unobservable with respect to the defined output. Hence, omitting $\phi_1'$ does not affect the output dynamics. Thus, we define $\phi''=[\phi_2', \dots, \phi_n']$, and obtain the dynamics
\begin{IEEEeqnarray}{lcl}
\dot{\hat{\omega}} &=& M \Big(- (K^{\text{droop}} + K^\omega) \hat{\omega} + K^V \hat{V} \nonumber \\
&& - {K^V}(K^\omega)^{-1} K^\text{droop, I} \eta - \mathcal{L}_\phi S \phi''  + P^{{m}}  \Big) \nonumber \\
\dot{\hat{V}} &=& 
\frac{1}{V^{\text{nom}}}E {K}^\omega \hat{\omega} -E\left(\mathcal{L}_R + \frac{K^V}{V^{\text{nom}}} \right) \hat{V} + \frac{1}{V^\text{nom}} E \mathcal{L}_\phi  S\phi''  \nonumber \\
\dot{\eta} &=&  K^{\text{droop,I}}\hat{\omega} - \mathcal{L}_\eta \eta \nonumber \\
\dot{\phi}'' &=&  S^T(K^V)^{-1}{K}^\omega \hat{\omega}  -\gamma \phi''. \label{eq:cl_dynamics_vec_delta_int_coordinated_double_prime}
\end{IEEEeqnarray}
 Provided that the system matrix, denoted $A$, of \eqref{eq:cl_dynamics_vec_delta_int_coordinated_double_prime} is full-rank, then \eqref{eq:cl_dynamics_vec_delta_int_coordinated_double_prime} has a unique equilibrium. Denote this equilibrium $x_0=[\omega_0^T, V_0^T, \eta_0^T, \phi_0''^T]^T$. Define $\bar{x}\triangleq [\bar{\omega}^T, \bar{V}^T, \bar{\eta}^T, \bar{\phi}^T]^T =[\hat{\omega}^T, \hat{V}^T, \eta^T, \phi''^T]^T - [\omega_0^T, V_0^T, \eta_0^T, \phi_0''^T]^T$. Hence, 
\begin{align}
\dot{\bar{x}} = A \bar{x},
\label{eq:dynamics_A_coordinated_shifted}
\end{align}
with the origin as its unique equilibrium. We are now ready to show the main stability result of this section. We first make the following assumptions.
\begin{assumption}
\label{ass:L_phi_coordinated}
The Laplacian matrix satisfies $\mathcal{L}_\phi = k_\phi \mathcal{L}_R$.
\end{assumption}
Assumption \ref{ass:L_phi_coordinated} can be interpreted as the emulated AC dynamics of \eqref{eq:voltage_control_secondary} having the same susceptance ratios as the conductance ratios of the HVDC lines. 
\begin{assumption}
\label{ass:gamma_coordinated}
The gain $\gamma$ satisfies $\gamma > \frac{k_\phi}{4V^\text{nom}}$. 
\end{assumption}
Assumption \ref{ass:gamma_coordinated} lower bounds for the damping coefficient $\gamma$.
\begin{theorem}
\label{th:stability_A_coordinated}
If Assumptions \ref{ass:L_phi_coordinated} and \ref{ass:gamma_coordinated} hold, the origin of \eqref{eq:dynamics_A_coordinated_shifted} is globally asymptotically stable. 
\end{theorem}
\begin{pf}
Consider the Lyapynov function candidate:
\begin{align}
W(\bar{\omega}, \bar{V}, \bar{\eta}, \bar{\phi}) &= \frac 12 \bar{\omega}^T K^\omega (K^V)^{-1} M^{-1}\bar{\omega} + \frac{V^\text{nom}}{2} \bar{V}^T C \bar{V} \nonumber \\
&\;\;\;\; + \frac 12 \bar{\eta}^T \bar{\eta} + \frac 12 \bar{\phi}^T S^T \mathcal{L}_\phi S \bar{\phi}, \label{eq:lyap_hvdc_coordinated}
\end{align}
where  $C=\diag(C_1, \dots, C_n)$. 
Clearly $W(\bar{\omega}, \bar{V}, \bar{\eta}, \bar{\phi})$ is positive definite and radially unbounded. Differentiating \eqref{eq:lyap_hvdc_coordinated} along trajectories of \eqref{eq:dynamics_A_coordinated_shifted}, we obtain
\begin{align*}
\dot{W}(\bar{\omega}, \bar{V}, \bar{\eta}, \bar{\phi})  &= \bar{\omega}^T K^\omega (K^V)^{-1} M^{-1}\dot{\bar{\omega}} + V^\text{nom} \bar{V}^T C \dot{\bar{V}} \nonumber \\
&\;\;\;\; + \bar{\eta}^T \dot{\bar{\eta}} + \bar{\phi}^T S^T \mathcal{L}_\phi S \dot{\bar{\phi}} \\
&= \bar{\omega}^T \big( -K^\omega (K^V)^{-1}(K^\omega + K^\text{droop})\bar{\omega} \\
&\;\;\;\; + K^\omega \bar{V} - K^\text{droop, I} \bar{\eta} - K^\omega (K^V)^{-1} \mathcal{L}_\phi S \bar{\phi} \big) \\
&\;\;\;\; + \bar{V}^T \Big( K^\omega \bar{\omega} - (V^\text{nom}\mathcal{L}_R {+} K^V)\bar{V} + \mathcal{L}_\phi S \bar{\phi} \Big) \\
&\;\;\;\; + \bar{\eta}^T \Big( K^\text{droop, I} \bar{\omega} - \mathcal{L}_\eta \eta  \Big) \\
&\;\;\;\; + \bar{\phi}^T S^T \mathcal{L}_\phi S \Big( S^T (K^V)^{-1}K^\omega \bar{\omega} - \gamma \bar{\phi} \Big) \\
&= -\bar{\omega}^T \big( -K^\omega (K^V)^{-1}(K^\omega + K^\text{droop})\bar{\omega} \\
&\;\;\;\; +  2 \bar{\omega}^T K^\omega \bar{V} - \bar{V}^T (V^\text{nom}\mathcal{L}_R + K^V)\bar{V} \\
&\;\;\;\; + \bar{V}^T \mathcal{L}_\phi S \bar{\phi} - \bar{\eta}^T \mathcal{L}_\eta \bar{\eta} - \gamma \bar{\phi}^T S^T \mathcal{L}_\phi S \bar{\phi}, 
\end{align*} 
since $SS^T = I_n - \frac 1n 1_{n\times n}$, so $\mathcal{L}_\phi SS^T = \mathcal{L}_\phi$. By defining 
\begin{align*}
\bar{V}' = \begin{bmatrix}
\frac{1}{\sqrt{n}} 1_n^T \\ S^T
\end{bmatrix} & \bar{V}
\qquad
\bar{V} = \begin{bmatrix}
\frac{1}{\sqrt{n}} 1_n & S
\end{bmatrix} \bar{V}',
\end{align*} 
we obtain $\bar{V}^T\mathcal{L}_R \bar{V} = \bar{V}''^T S^T\mathcal{L}_R S \bar{V}'' $ and $\bar{V}^T\mathcal{L}_\phi S \bar{\phi} = \bar{V}''^T S^T\mathcal{L}_\phi S \bar{\phi}$, where $\bar{V}'' = [\bar{V}'_2, \dots, \bar{V}'_n]^T$. By invoking Assumption~\ref{ass:L_phi_coordinated}, we obtain
\begin{align*}
&\dot{W}(\bar{\omega}, \bar{V}, \bar{\eta}, \bar{\phi})  \\
&= - \begin{bmatrix}
\bar{\omega}^T & \bar{V}^T 
\end{bmatrix}
\underbrace{\begin{bmatrix}
K^\omega (K^V)^{-1}(K^\omega + K^\text{droop}) & -K^\omega \\
-K^\omega & K^V 
\end{bmatrix}}_{\triangleq Q_1}
\begin{bmatrix}
\bar{\omega} \\ \bar{V}
\end{bmatrix} \\
&\;\;\;\; - \begin{bmatrix}
\bar{V}''^T & \bar{\phi}^T
\end{bmatrix}
\underbrace{\begin{bmatrix}
V^\text{nom} S^T\mathcal{L}_R S & -\frac{k_\phi}{2} S^T\mathcal{L}_R S \\
-\frac{k_\phi}{2} S^T\mathcal{L}_R S & \gamma k_\phi S^T\mathcal{L}_R S
\end{bmatrix}}_{\triangleq Q_2}
\begin{bmatrix}
\bar{V}'' \\ \bar{\phi}
\end{bmatrix} \\
&\;\;\;\; - \bar{\eta}^T \mathcal{L}_\eta \bar{\eta}.  
\end{align*}
By applying the Schur complement condition for positive definiteness, we see that $Q_1$ is positive definite iff
\begin{eqnarray*}
K^\omega (K^V)^{-1}(K^\omega + K^\text{droop}) - K^\omega (K_V)^{-1} K^\omega \\
= K^\omega (K^V)^{-1} K^\text{droop} > 0.
\end{eqnarray*}
Hence $Q_1$ is positive definite. By applying the same argument to $Q_2$, $Q_2$ is positive definite iff
\begin{align*}
\Bigg( \gamma k_\phi- \frac{k_\phi^2}{4 V^\text{nom}} \Bigg) S^T \mathcal{L}_R S > 0.
\end{align*}
Clearly the above matrix inequality holds under Assumption~\ref{ass:gamma_coordinated}, since $S^T \mathcal{L}_R S \ge 0$, and $Sx\ne k 1_n$ for $k\ne 0$. Thus $\dot{W}(\bar{\omega}, \bar{V}, \bar{\eta}, \bar{\phi}) \le 0$, and the set where $W(\bar{\omega}, \bar{V}, \bar{\eta}, \bar{\phi})$ is non-decreasing is given by
\begin{align*}
G &= \{ (\bar{\omega}, \bar{V}, \bar{\eta}, \bar{\phi}) | \dot{W}(\bar{\omega}, \bar{V}, \bar{\eta}, \bar{\phi}) =0 \} \\ 
&=  \{ (\bar{\omega}, \bar{V}, \bar{\eta}, \bar{\phi}) | \bar{\eta} = k 1_n \},
\end{align*}
for any $k\in \mathbb{R}$. Clearly the largest invariant set in $G$ is the origin, and thus $k=0$. By LaSalle's theorem, the origin is globally asymptotically stable under the dynamics \eqref{eq:dynamics_A_coordinated_shifted}. 
\end{pf}
We now turn our attention to the equilibrium of  \eqref{eq:cl_dynamics_vec_delta_int_coordinated_double_prime}. If the stability condition in Theorem~\ref{th:stability_A_coordinated} is met, then the equilibrium  is stable, and $A$ is thus Hurwitz. Hence \eqref{eq:cl_dynamics_vec_delta_int_coordinated_double_prime} must have a unique equilibrium. In the following corollary we show that the equilibrium of \eqref{eq:cl_dynamics_vec_delta_int_coordinated_double_prime} satisfies Objective~\ref{obj:1_hvdc_coordinated}.
\begin{corollary}
\label{cor:hvdc_coordinated_equilibrium}
Let Assumption~\ref{ass:L_phi_coordinated} hold and let $\gamma$, $k_\phi$ be given such that Assumption~\ref{ass:gamma_coordinated} holds. 
Let $K^V, K^\omega$ and $K^\text{droop}$ be such that $(F^P)^{-1} = K^V(K^\omega)^{-1}K^\text{droop}$ and $F^V=K^V$, where $F^P = \diag(f^P_1, \dots, f^P_n)$ and $F^V = \diag(f^V_1, \dots, f^V_n)$. 
 Then Objective~\ref{obj:1_hvdc_coordinated} is satisfied in the limit when $\norm{(K^\omega)^{-1} K^V}_\infty \rightarrow 0$. 
\end{corollary}
\begin{pf}
By Theorem~\ref{th:stability_A_coordinated},  \eqref{eq:cl_dynamics_vec_delta_int_coordinated_double_prime} has a unique and  stable equilibrium. 
The last $n-1$ rows of the equilibrium imply
\begin{align*}
S^T(K^V)^{-1}K^\omega \hat{\omega} - \gamma \phi'' = 0_{n-1}.
\end{align*}
Now $\norm{(K^\omega)^{-1} K^V}_\infty \rightarrow 0$ in the above equation implies that $S^T\hat{\omega} = 0 \Leftrightarrow \hat{\omega} = k_1 1_n$ for some $k_1\in \mathbb{R}$. Consider the $(2n+1)$th to $3n$th rows of the equilibrium of \eqref{eq:cl_dynamics_vec_delta_int_coordinated_double_prime}:
\begin{align*}
K^\text{droop, I} \hat{\omega} - \mathcal{L}_\eta \eta = 0_n.
\end{align*}
By inserting $\hat{\omega} = k_1 1_n$ and premultiplying the above equation with $1_n^T$, we obtain that $k_1=0$, so $\hat{\omega}=0_n$ and Equation~\eqref{eq:hvdc_coordinated_frequency_objective} of Objective~\ref{obj:1_hvdc_coordinated} is satisfied. This implies that $\eta = k_2 1_n$ for some $k_2\in \mathbb{R}$.  Finally we consider the $(n+1)$th to $2n$th rows of the equilibrium of \eqref{eq:cl_dynamics_vec_delta_int_coordinated_double_prime}:
\begin{align}
\frac{1}{V^\text{nom}} K^\omega \hat{\omega} - \Big( \mathcal{L}_R + \frac{K^V}{V^\text{nom}} \Big) \hat{V} + \frac{1}{V^\text{nom}} \mathcal{L}_\phi S \phi'' = 0_n \label{eq:cl_dynamics_vec_delta_int_coordinated_double_prime_n+1:2n}
\end{align}
Inserting $\hat{\omega}=0_n$ and premultiplying \eqref{eq:cl_dynamics_vec_delta_int_coordinated_double_prime_n+1:2n} with $1_n^T$ yield  
\begin{align}
1_n^TK^V\hat{V}=0.
\label{eq:hvdc_coordinated_voltage_equilibrium}
\end{align}
Inserting $\hat{\omega}=0_n$ and $\eta = k_2 1_n$ in \eqref{eq:hvdc_coordinated_droop_control_secondary_distributed} yields 
\begin{align}
P^\text{gen} = - k_2 K^V(K^\omega)^{-1} K^\text{droop, I} 1_n,
\label{eq:hvdc_coordinated_eta_equilibrium}
\end{align}
where $P^\text{gen} = [P^\text{gen}_1, \dots, P^\text{gen}_n]^T$.
It now remains to show that the equilibrium of \eqref{eq:cl_dynamics_vec_delta_int_coordinated_double_prime} minimizes the  cost functions \eqref{eq:hvdc_coordinated_generation_objective}  and \eqref{eq:hvdc_coordinated_voltage_objective}. Consider first \eqref{eq:hvdc_coordinated_generation_objective}, with the constraints
$P^\text{gen}_i + P_i^{{m}} - P_i^{\text{inj}} = 0 \; \forall i=1, \dots, n$ and $\sum_{i=1}^n P_i^\text{inj} = 0$. By summing the first constraints we obtain $\sum_{i=1}^n P^\text{gen}_i  = - \sum_{i=1}^n P_i^{{m}} $. The KKT condition of \eqref{eq:hvdc_coordinated_generation_objective} is 
\begin{align}
F^P P^\text{gen} = - k_3 1_n.
\label{eq:hvdc_coordinated_KKT_generation}
\end{align}
Since $(F^P)^{-1} = K^V(K^\omega)^{-1}K^\text{droop}$,  \eqref{eq:hvdc_coordinated_eta_equilibrium} and \eqref{eq:hvdc_coordinated_KKT_generation} are identical for $k_2=k_3$. We conclude that  \eqref{eq:hvdc_coordinated_generation_objective} is minimized. 
Since $P^\text{gen} =  - K^V(K^\omega)^{-1} K^\text{droop} \eta = - k_2 K^V(K^\omega)^{-1} K^\text{droop} 1_n$ and $\hat{\omega}=0_n$, premultiplying the first $n$ rows of the equilibrium of \eqref{eq:cl_dynamics_vec_delta_int_coordinated_double_prime} with $M^{-1}$, and adding to the $(n+1)$th to $2n$th rows premultiplied with $V^\text{nom}E^{-1}$ yields
\begin{align*}
- V^\text{nom} \mathcal{L}_R \hat{V} - k_2 K^V(K^\omega)^{-1} K^\text{droop} 1_n = P^m.
\end{align*}
Premultiplying the above equation with $1^T_n$ yields $k_2 = - \sum_{i=1}^n P^m_i \sum_{i=1}^n \frac{K^\omega_i}{K^V_i K^\text{droop}_i}$. Additionally, $\mathcal{L}_R \hat{V}$ is uniquely determined. 
Now consider \eqref{eq:hvdc_coordinated_voltage_objective}. Note that $P_i^\text{inj}$ and hence $I_i^\text{inj}$, are uniquely determined by \eqref{eq:hvdc_coordinated_generation_objective}. By the equilibrium of \eqref{eq:hvdc_coordinated_voltage}, $\mathcal{L}_R \hat{V} = I^\text{inj}$, where $I^\text{inj} = [I^\text{inj}_1, \dots, I^\text{inj}_n]^T$. Thus, the KKT condition of \eqref{eq:hvdc_coordinated_voltage_objective} is 
\begin{align}
F^V \hat{V} =  \mathcal{L}_R r,
\label{eq:hvdc_coordinated_KKT_voltage_raw}
\end{align}
where $r\in \mathbb{R}^n$. 
Since $\mathcal{L}_R \hat{V}$ is uniquely determined, we premultiply \eqref{eq:hvdc_coordinated_KKT_voltage_raw} with $1_n^T$ and obtain the equivalent condition
\begin{align}
1^T_n F^V \hat{V} = 0.
\label{eq:hvdc_coordinated_KKT_voltage}
\end{align}
Since $F^V = K^V$, \eqref{eq:hvdc_coordinated_voltage_equilibrium} and \eqref{eq:hvdc_coordinated_KKT_voltage} are equivalent. Hence \eqref{eq:hvdc_coordinated_voltage_objective} is minimized, so Objective~\ref{obj:1_hvdc_coordinated} is satisfied. 
\end{pf}

\begin{remark}
Corollary \ref{cor:hvdc_coordinated_equilibrium} provides insight in choosing the controller gains of \eqref{eq:hvdc_coordinated_droop_control_secondary_distributed} and \eqref{eq:voltage_control_secondary}, to satisfy Objective~\ref{obj:1_hvdc_coordinated}. 
\end{remark}

\section{Simulations}
\label{sec:simulations}
\begin{figure}[ht]
\def\svgwidth{1.25\columnwidth}
\begin{center}
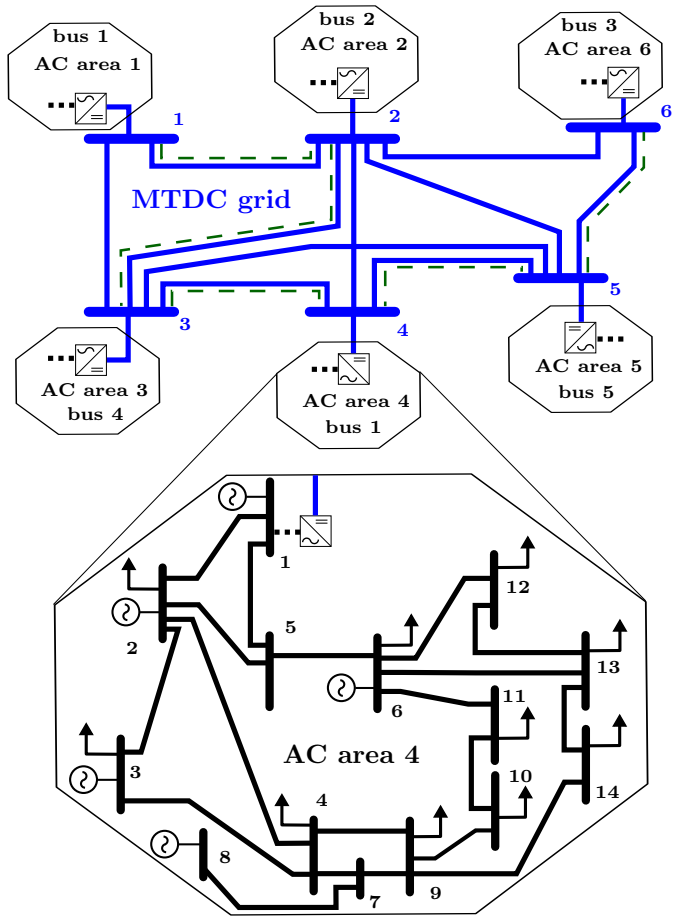
\captionsetup{width=\linewidth}	
	\caption{MTDC test system, consisting of a 6-terminal MTDC grid. Each terminal is connected to an IEEE 14 bus AC grid, sketched as octagons.}
	\label{fig:testsystem}
\end{center} 
\end{figure}

In this section, simulations are conducted on a test system to validate the performance of the proposed controllers. The simulation was performed in Matlab, using a dynamic phasor approach based on \citep{demiray2008}. The test system is illustrated in Figure~\ref{fig:testsystem}.
The line parameters of the MTDC grid are given in Table~\ref{tab:HVDCgridParameter}. Note that we in the simulation also consider the inductances $L_{ij}$ and capacitances $C_{ij}$ of the HVDC lines. The capacitances of the terminals are assumed to be given by $C_i=0.375\times 10^{-3}$ p.u., and  are chosen uniformly for all VSC stations. The AC grid parameters were obtained from \citep{Milano}. The generators are modeled as a 6{th} order machine model controlled by an automatic voltage controller and a governor \citep{kundur1994power}. Each load in the grid is assumed to be equipped with an ideal power controller.

\begin{table}[htb]
\captionsetup{width=\linewidth}
\caption{HVDC grid line parameters}
\label{tab:HVDCgridParameter}	
\begin{center}
\begin{tabular}{ccccc} 
		$i$&$j$&$R_{ij}$ [p.u.]&$L_{ij}$ [$10^{-3}$ p.u.]& $C_{ij}$ [p.u.]\\ \hline
		1 & 2 & 0.0586 & 0.2560 &0.0085 \\
		1 & 3 & 0.0586 & 0.2560 &0.0085\\
		2 & 3 & 0.0878 & 0.3840 &0.0127\\
		2 & 4 & 0.0586 & 0.2560 &0.0085\\
		2 & 5 & 0.0732 & 0.3200 &0.0106\\
		2 & 6 & 0.1464 & 0.6400 &0.0212\\
		3 & 4 & 0.0586 & 0.2560 &0.0085\\
		3 & 5 & 0.1464 & 0.6400 &0.0212\\
		4 & 5 & 0.0732 & 0.3200 &0.0106\\
		5 & 6 & 0.1464 & 0.6400 &0.0212
\end{tabular}
\end{center}
\end{table}

\begin{table}[htb]
\begin{center}			
\captionsetup{width=\linewidth}	
\caption{Controller parameters}
\label{tab:ControllerParameter}
\begin{tabular}{ccccccc}
$K^{\omega}_i$ &$K^{V}_i$&$K^\text{droop}_i$& $K^\text{droop, I}_i$ &$c^{\eta}_{ij}$ &$c^{\phi}_{ij}$ & $\gamma$ \\
\hline
1000            &   100       & 9   & 3.35      &  ${5}/{R_{ij}}$ & ${15}/{R_{ij}}$  &  0                               
\end{tabular}
\end{center}
\end{table}

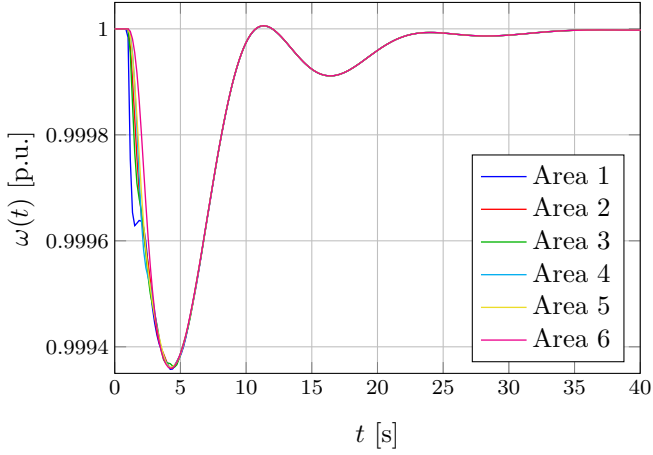
\begin{figure}[t!]
	\centering
\tikzsetnextfilename{Frequencies}
	\begin{tikzpicture}
	\begin{axis}
	[cycle list name=Necsys,
	xlabel={$t$ [s]},
	ylabel={$\omega(t)$ [p.u.]},
	xmin=0,
	xmax=40,
	xtick={0,5,...,40},
	ymin=0.99935,
	ymax=1.00005,
	yticklabel style={/pgf/number format/.cd,
		fixed,
		precision=4},
	grid=major,
	height=6.5cm,
	width=0.96\columnwidth,
	legend cell align=left,
	legend pos= south east,
	]
	\foreach \x in {1,2,3,4,5,6}{
	\addplot table[x index=0,y index=\x,col sep=tab]{PlotData_coordinated/Necsys_Frequency.txt};
	}	
	\legend{Area 1, Area 2, Area 3, Area 4, Area 5, Area 6}
	\end{axis}
	\end{tikzpicture}
	\caption{Average frequencies in the AC areas.}
	\label{fig:ACFrequency}
\end{figure}

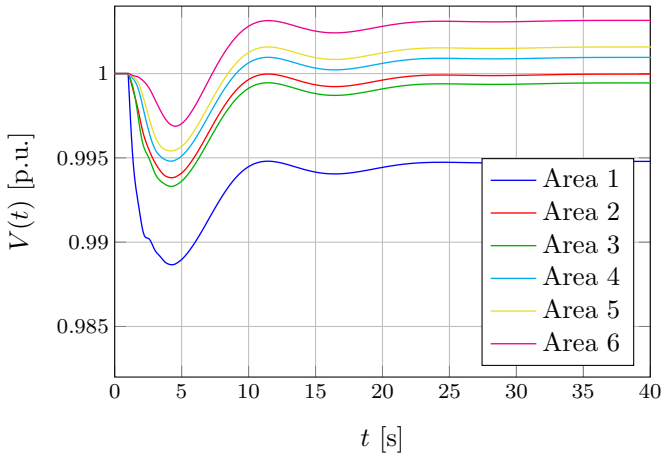
\begin{figure}[htb]
	\centering
\tikzsetnextfilename{DCVoltages}
	\begin{tikzpicture}
	\begin{axis}[cycle list name=Necsys,
	xlabel={$t$ [s]},
	ylabel={$V(t)$ [p.u.]},
	xmin=0,
	xmax=40,
	xtick={0,5,...,40},
	ymin=0.982,
	ymax=1.004,
	yticklabel style={/pgf/number format/.cd,
		fixed,
		precision=4},
	grid=major,
	height=6.5cm,
	width=0.975\columnwidth,
	legend pos= south east,
	]
	\foreach \x in {1,2,3,4,5,6}{
	\addplot table[x index=0,y index=\x,col sep=tab] {PlotData_coordinated/Necsys2_Voltage.txt};
	}
	\legend{Area 1, Area 2, Area 3, Area 4, Area 5, Area 6}
	\end{axis}
	\end{tikzpicture}
	\caption{DC terminal voltages.}
	\label{fig:DCVoltage}
\end{figure}

\begin{figure}[t!]
	\centering
\tikzsetnextfilename{Generators}
	\begin{tikzpicture}
	\begin{axis}[cycle list name=Necsys,
	xlabel={$t$ [s]},
	ylabel={$P^\text{gen}(t)$ [p.u.]},
	xmin=0,
	xmax=40,
	xtick={0,5,...,40},
	ymin=-0.005,
	ymax=0.076,
	scaled ticks=false,
	tick label style={/pgf/number format/fixed,
		/pgf/number format/precision=2},
	ytick={0,0.02,...,0.101},
	grid=major,
	height=6.5cm,
	width=0.975\columnwidth,
	legend pos= north east,
	]
	\foreach \x in {1,2,3,4,5,6}{
		\addplot table[x index=0,y index=\x,col sep=tab] {PlotData_coordinated/Necsys2_Generators.txt};
	}
	\legend{Area 1, Area 2, Area 3, Area 4, Area 5, Area 6}
	\end{axis}
	\end{tikzpicture} 
	\caption{Increase of generated power in the AC areas.}
	\label{fig:Generators}
\end{figure}
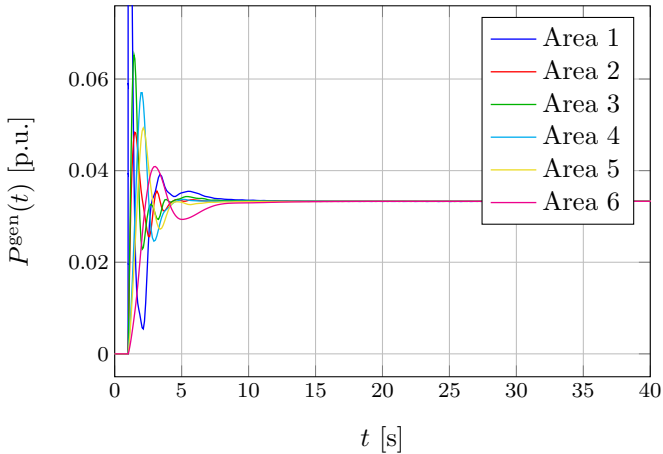

The test grid was controlled with the controllers \eqref{eq:hvdc_coordinated_droop_control_secondary_distributed} and \eqref{eq:voltage_control_secondary}, with parameters given in Table~\ref{tab:ControllerParameter}.
The communication network of  \eqref{eq:hvdc_coordinated_droop_control_secondary_distributed} and \eqref{eq:voltage_control_secondary} is illustrated by the dashed lines in Figure~\ref{fig:testsystem}. 
 Note that we have set $\gamma=0$, and that Theorem~\ref{th:stability_A_coordinated} thus does not guarantee stability of the equilibrium. The resulting matrix of the closed-loop system is however verified to be Hurwitz. At time $t=1$ the output of one generator in area 1 was reduced by $0.2$ p.u., simulating a fault. 
Figure~\ref{fig:ACFrequency} shows the average frequencies of the AC grids, while Figure~\ref{fig:DCVoltage} shows the DC voltages of the terminals. Figure~\ref{fig:Generators} shows the total increase of the generated power within each AC area.
It can be noted that immediately after the fault, the average frequency of the AC area of the fault drops. The frequency drop is followed by a voltage drop in all terminals, and a subsequent frequency drop in the remaining AC areas. The frequencies converge to the nominal frequency, while the voltages converge to their new stationary values after approximately $30$ s. We note that the frequencies are restored to the nominal frequency, as predicted by Corollary~\ref{cor:hvdc_coordinated_equilibrium}. Furthermore the incremental generated power of the AC areas converge to the same value, as a consequence of Corollary~\ref{cor:hvdc_coordinated_equilibrium} and the controller parameters being equal.

\section{Discussion and Conclusions}
\label{sec:discussion}
In this paper we have studied a distributed secondary controller for sharing frequency control reserves of asynchronous AC systems connected through an MTDC system. Under certain conditions, the proposed controller stabilizes the closed-loop system consisting of the interconnected AC systems and the MTDC grid. The frequencies in all AC grids are shown to converge to the nominal frequency. Furthermore, quadratic cost functions of the voltage deviations of the MTDC terminals and of the generated power, are minimized asymptotically. Finally, the results were validated on a six-terminal MTDC system with connected AC systems.

\bibliography{references}
\bibliographystyle{plain}
\end{document}

%% file: ACandDCGrid2.pdf_tex
\begingroup%
  \makeatletter%
  \providecommand\color[2][]{%
    \errmessage{(Inkscape) Color is used for the text in Inkscape, but the package 'color.sty' is not loaded}%
    \renewcommand\color[2][]{}%
  }%
  \providecommand\transparent[1]{%
    \errmessage{(Inkscape) Transparency is used (non-zero) for the text in Inkscape, but the package 'transparent.sty' is not loaded}%
    \renewcommand\transparent[1]{}%
  }%
  \providecommand\rotatebox[2]{#2}%
  \ifx\svgwidth\undefined%
    \setlength{\unitlength}{649.68979373bp}%
    \ifx\svgscale\undefined%
      \relax%
    \else%
      \setlength{\unitlength}{\unitlength * \real{\svgscale}}%
    \fi%
  \else%
    \setlength{\unitlength}{\svgwidth}%
  \fi%
  \global\let\svgwidth\undefined%
  \global\let\svgscale\undefined%
  \makeatother%
  \begin{picture}(1,1.09103572)%
    \put(0,0){\includegraphics[width=\unitlength,page=1]{ACandDCGrid2.pdf}}%
    \put(0.19686414,0.94820045){\color[rgb]{0,0,0}\makebox(0,0)[lb]{\smash{\textbf{\textcolor{blue}{\scriptsize{1}}}}}}%
    \put(0,0){\includegraphics[width=\unitlength,page=2]{ACandDCGrid2.pdf}}%
    \put(0.0314905,1.01318109){\color[rgb]{0,0,0}\makebox(0,0)[lb]{\smash{\textbf{\scriptsize{AC area 1}}}}}%
    \put(-0.0106717,1.19328518){\color[rgb]{0,0,0}\makebox(0,0)[lt]{\begin{minipage}{0.16636605\unitlength}\raggedright \end{minipage}}}%
    \put(0,0){\includegraphics[width=\unitlength,page=3]{ACandDCGrid2.pdf}}%
    \put(0.45371901,0.95066256){\color[rgb]{0,0,0}\makebox(0,0)[lb]{\smash{\textbf{\textcolor{blue}{\scriptsize{2}}}}}}%
    \put(0.7767853,0.95774878){\color[rgb]{0,0,0}\makebox(0,0)[lb]{\smash{\textbf{\textcolor{blue}{\scriptsize{6}}}}}}%
    \put(0.20407638,0.70193715){\color[rgb]{0,0,0}\makebox(0,0)[lb]{\smash{\textbf{\textcolor{blue}{\scriptsize{3}}}}}}%
    \put(0.46311239,0.69476366){\color[rgb]{0,0,0}\makebox(0,0)[lb]{\smash{\textbf{\textcolor{blue}{\scriptsize{4}}}}}}%
    \put(0.7160707,0.74338959){\color[rgb]{0,0,0}\makebox(0,0)[lb]{\smash{\textbf{\textcolor{blue}{\scriptsize{5}}}}}}%
    \put(0,0){\includegraphics[width=\unitlength,page=4]{ACandDCGrid2.pdf}}%
    \put(0.03946648,0.62260061){\color[rgb]{0,0,0}\makebox(0,0)[lb]{\smash{\textbf{\scriptsize{AC area 3}}}}}%
    \put(0,0){\includegraphics[width=\unitlength,page=5]{ACandDCGrid2.pdf}}%
    \put(0.34836246,1.03802671){\color[rgb]{0,0,0}\makebox(0,0)[lb]{\smash{\textbf{\scriptsize{AC area 2}}}}}%
    \put(0,0){\includegraphics[width=\unitlength,page=6]{ACandDCGrid2.pdf}}%
    \put(0.6273616,0.65214138){\color[rgb]{0,0,0}\makebox(0,0)[lb]{\smash{\textbf{\scriptsize{AC area 5}}}}}%
    \put(0,0){\includegraphics[width=\unitlength,page=7]{ACandDCGrid2.pdf}}%
    \put(0.63812943,1.0317465){\color[rgb]{0,0,0}\makebox(0,0)[lb]{\smash{\textbf{\scriptsize{AC area 6}}}}}%
    \put(0.76456297,0.80924392){\color[rgb]{0,0,0}\makebox(0,0)[lb]{\smash{}}}%
    \put(0,0){\includegraphics[width=\unitlength,page=8]{ACandDCGrid2.pdf}}%
    \put(0.3235102,0.4169686){\color[rgb]{0,0,0}\makebox(0,0)[lb]{\smash{\textbf{\scriptsize{1}}}}}%
    \put(0.14123496,0.31384403){\color[rgb]{0,0,0}\makebox(0,0)[lb]{\smash{\textbf{\scriptsize{2}}}}}%
    \put(0.14534406,0.16535265){\color[rgb]{0,0,0}\makebox(0,0)[lb]{\smash{\textbf{\scriptsize{3}}}}}%
    \put(0.4299907,0.00898873){\color[rgb]{0,0,0}\makebox(0,0)[lb]{\smash{\textbf{\scriptsize{7}}}}}%
    \put(0.33000951,0.33856181){\color[rgb]{0,0,0}\makebox(0,0)[lb]{\smash{\textbf{\scriptsize{5}}}}}%
    \put(0.36672831,0.13103624){\color[rgb]{0,0,0}\makebox(0,0)[lb]{\smash{\textbf{\scriptsize{4}}}}}%
    \put(0.45626113,0.23813242){\color[rgb]{0,0,0}\makebox(0,0)[lb]{\smash{\textbf{\scriptsize{6}}}}}%
    \put(0.59296117,0.3867491){\color[rgb]{0,0,0}\makebox(0,0)[lb]{\smash{\textbf{\scriptsize{12}}}}}%
    \put(0.70050551,0.29022504){\color[rgb]{0,0,0}\makebox(0,0)[lb]{\smash{\textbf{\scriptsize{13}}}}}%
    \put(0.69923041,0.14099863){\color[rgb]{0,0,0}\makebox(0,0)[lb]{\smash{\textbf{\scriptsize{14}}}}}%
    \put(0.58853673,0.25641714){\color[rgb]{0,0,0}\makebox(0,0)[lb]{\smash{\textbf{\scriptsize{11}}}}}%
    \put(0.59587204,0.15961997){\color[rgb]{0,0,0}\makebox(0,0)[lb]{\smash{\textbf{\scriptsize{10}}}}}%
    \put(0.50255767,0.02045717){\color[rgb]{0,0,0}\makebox(0,0)[lb]{\smash{\textbf{\scriptsize{9}}}}}%
    \put(0.25164642,0.06297213){\color[rgb]{0,0,0}\makebox(0,0)[lb]{\smash{\textbf{\scriptsize{8}}}}}%
    \put(0,0){\includegraphics[width=\unitlength,page=9]{ACandDCGrid2.pdf}}%
    \put(0.14743732,0.84840804){\color[rgb]{0,0,0}\makebox(0,0)[lb]{\smash{\textbf{\textcolor{blue}{MTDC grid}}}}}%
    \put(0,0){\includegraphics[width=\unitlength,page=10]{ACandDCGrid2.pdf}}%
    \put(0.32772548,0.18241149){\color[rgb]{0,0,0}\makebox(0,0)[lb]{\smash{\textbf{AC area 4}}}}%
    \put(0,0){\includegraphics[width=\unitlength,page=11]{ACandDCGrid2.pdf}}%
    \put(0.07092236,0.59451006){\color[rgb]{0,0,0}\makebox(0,0)[lb]{\smash{\textbf{\scriptsize{bus 4}}}}}%
    \put(0,0){\includegraphics[width=\unitlength,page=12]{ACandDCGrid2.pdf}}%
    \put(0.05370394,1.04673664){\color[rgb]{0,0,0}\makebox(0,0)[lb]{\smash{\textbf{\scriptsize{bus 1}}}}}%
    \put(0.37020368,1.0671981){\color[rgb]{0,0,0}\makebox(0,0)[lb]{\smash{\textbf{\scriptsize{bus 2}}}}}%
    \put(0.65840564,1.06023249){\color[rgb]{0,0,0}\makebox(0,0)[lb]{\smash{\textbf{\scriptsize{bus 3}}}}}%
    \put(0.6548503,0.62096394){\color[rgb]{0,0,0}\makebox(0,0)[lb]{\smash{\textbf{\scriptsize{bus 5}}}}}%
    \put(0.37804001,0.57829965){\color[rgb]{0,0,0}\makebox(0,0)[lb]{\smash{\textbf{\scriptsize{bus 1}}}}}%
    \put(0,0){\includegraphics[width=\unitlength,page=13]{ACandDCGrid2.pdf}}%
    \put(0.35070136,0.60964889){\color[rgb]{0,0,0}\makebox(0,0)[lb]{\smash{\textbf{\scriptsize{AC area 4}}}}}%
  \end{picture}%
\endgroup%